\documentclass[a4paper,11pt,leqno]{article}

\usepackage{amssymb,latexsym}
\usepackage{mltex}
\usepackage{amsmath}
\usepackage{fancyhdr}
\usepackage{enumerate}
\usepackage{amsthm}
\usepackage{hyperref}

\newtheorem{thm}{Theorem}     
\newtheorem{cor}{Corollary}

\newtheorem{rem}{Remark}

\newtheorem{thrm}{Theorem}

\newcommand{\vvol}{\mathrm{Vol}}

\newcommand{\nablab}{\overline{\nabla}}

\newcommand{\Ric}{\mathrm{Ric}}

\newcommand{\lgra}{\longrightarrow}
\newcommand{\Hinf}{||H||_{\infty}}

\newcommand{\iid}{\mathrm{Id}\,}

\newcommand{\sscal}{\mathrm{Scal}\,}

\newcommand{\trace}{\mathrm{tr\,}}

\newcommand{\Ss}{\mathbb{S}}

\newcommand{\R}{\mathbb{R}}

\newcommand{\beqt}{\begin{equation}}  \newcommand{\eeqt}{\end{equation}}
\newcommand{\bal}{\begin{align}}      \newcommand{\eal}{\end{align}}
\newcommand{\ba}{\begin{array}}      \newcommand{\ea}{\end{array}}
\newcommand{\bc}{\begin{center}}     \newcommand{\ec}{\end{center}}
\newcommand{\be}{\begin{enumerate}}  \newcommand{\ee}{\end{enumerate}}
\newcommand{\beq}{\begin{eqnarray}}  \newcommand{\eeq}{\end{eqnarray}}
\newcommand{\beQ}{\begin{eqnarray*}} \newcommand{\eeQ}{\end{eqnarray*}}
\newcommand{\bi}{\begin{itemize}}    \newcommand{\ei}{\end{itemize}}
\newcommand{\bt}{\begin{tabular}}    \newcommand{\et}{\end{tabular}}
\newcommand{\finpreuve}{\hfill\square\\}

\def\pf{\noindent{\textit {Proof :} }}

\title{A Remark on Almost Umbilical Hypersurfaces}
\author{Julien Roth}
\date{\today}
\begin{document}

\maketitle
\begin{center}
Universit\'e Paris-Est, LAMA (UMR 8050), UPEMLV, UPEC, CNRS, F-77454, Marne-la-Vall\'ee, France
\end{center}
\begin{center}
julien.roth@univ-mlv.fr
\end{center}
\begin{abstract}
In this article, we prove new stability results for almost-Einstein hypersurfaces of the Euclidean space, based on previous eigenvalue pinching results. Then, we deduce some comparable results for almost umbilical hypersurfaces. 
\end{abstract}
\noindent
{\bf Mathematical Subject Clasification:} 53A07, 53A10, 53C20, 53C24.

\date{}
\maketitle\pagenumbering{arabic}

\section{Introduction}\label{intro}
It is a well-known fact that a totally umbilical hypersurface of the Euclidean space $\R^{n+1}$ which is not totally geodesic is a round sphere. An other classical rigidty result states that an Einstein (with positive scalar curvature) hypersurface of the Euclidean space $\R^{n+1}$ is a round sphere. This was proved by Thomas \cite{Th} and independently by Fialkow \cite{Fi} in the 30's. Recently, Grosjean \cite{Gr} gave a new proof based on the equality case of an estimate of the first eigenvalue of the Laplacian involving the scalar curvature.\\
It is then natural to consider the case of almost umbilical and almost Einstein hypersurfaces and ask if they are close to round spheres (in a sense to be precised). Shiohama and Xu \cite{SX1,SX2} proved that under a condition on Betti numbers, almost umbilical hypersurfaces of Euclidean space are homeomorphic to the sphere. For almost Einstein hypersurfaces, Vlachos obtained a comparable result in \cite{Vl}. \\
By studying the stablity problem of Grosjean's proof, we showed in \cite{Ro} another promity result for almost Einstein, namely these hypersurfaces are diffeomorphic and quasi-isometric to round sphere. We got stronger proximity but in counterpart, with stronger assumption that in the result of Vlachos.\\
The aim of the present paper is to relax the assumptions of this result in one hand (Theorem \ref{thm1}) and on the other hand to get stability results for almost umbilical hypersurfaces to be compared with the results of Shiohama and Xu (see Theorems \ref{thm2} and \ref{thm3}).
\indent

\section{Preliminaries}\label{prelim}
Let $(M^n,g)$ be a $n$-dimensional compact, connected, oriented Riemannian manifold without boundary, isometrically immersed into the $(n+1)$-dimensional Euclidean space $(\R^{n+1},can)$. The second fundamental form $B$ of the immersion is the bilinear symmetric form defined by
$$B(Y,Z)=-g\left( \nablab_Y\nu,Z\right),$$
where $\nablab$ is the Riemannian connection on $\R^{n+1}$ and $\nu$ the outward normal unit vector field on $M$.\\
\indent
From $B$, we can define the mean curvature,
$$H=\frac{1}{n}\trace(B),$$
and, more generally, the higher order mean curvatures,
$$H_r=\frac{1}{\left( \begin{array}{c}
 n\\r 
\end{array}\right) }\sigma_r(\kappa_1,\cdots,\kappa_n),$$
where $\sigma_r$ is the $r$-th symmetric polynomial and $\kappa_1,\cdots,\kappa_n$ are the principal curvatures of the immersion. By convention, we set $H_0=1$.\\
\indent
Note that $H_1=H$ and from the Gauss equation, $H_2$ is, up to a mutliplicative constant, the scalar curvature. Namely, we have $H_2=\frac{1}{n(n-1)}\sscal$.\\
\indent
These extrinsic curvatures satisfy the well-known Hsiung-Minkowski formula, for $1\leqslant r\leqslant n$,
\beqt\label{Hsuing}
\int_M\big( H_{r-1}+H_r\left\langle X,\nu\right\rangle \big)dv_g=0,
\eeqt
where $X$ is the position vector and $\nu$ the normal vector of the immersion. They also satisfy the following inequalities if $H_r$ is a positive function:
\beqt\label{ineq}
H_r^{\frac{1}{r}}\leq H_{r-1}^{\frac{1}{r-1}}\leq\cdots\leq
H_2^{\frac{1}{2}}\leq H.
\eeqt

Moreover, Reilly \cite{Re} proved some upper bounds for the first eigenvalue of the Laplacian for hypersurfaces of $\R^{n+1}$ in terms of higher order mean curvatures. Precisely, he showed
\beqt\label{reilly2}
\lambda_1(M)\left( \int_MH_{r-1}dv_g\right)^2 \leqslant\frac{n}{\vvol(M)}\int_MH_r^2dv_g,
\eeqt 
with equality if and only if $m$ is a geodesic hyperspheres of $\R^{n+1}$.\\
\indent
By the H\"older inequality, we obtain for $p\geqslant2$, 
$$ \lambda_1(M)\leqslant n\frac{||H_r||_{2p}^2\vvol(M)^{2}}{\left( \int_MH_{r-1}dv_g\right)^2}.$$
Now, for $p\geqslant2$ and $1\leqslant r\leqslant n$, we define $k_{p,r}=\dfrac{||H_r||_{2p}^2\vvol(M)^{2}}{\left( \int_MH_{r-1}dv_g\right)^2}$, which are the constants involved in Theorem \ref{thm1}.\\
\indent
Note that we use the following convention, for any smooth fonction $f$ and any $p\geqslant1$,
$$||f||_p=\frac{1}{\vvol(M)^{\frac{1}{p}}}\left(\int_M|f|^pdv_g\right)^{\frac{1}{p}}.$$
The main tool that we will use in the present paper is the following pinching result, associated with these inequalities, that we proved in \cite{Ro}.
\begin{thrm}[Roth \cite{Ro}]\label{thrm2}
 Let $(M^n,g)$ be a compact, connected, oriented Riemannian manifold without boundary isometrically immersed in $\R^{n+1}$, $n\geqslant2$.  Let $r\in\{1,\cdots,n\}$ and assume that $H_r>0$ if $r>1$. Then for any $q\geqslant2$ and any $\theta\in]0,1[$, there exists a constant $K_{\theta}$ depending only on $n$, $\Hinf$, $\vvol(M)$, $||H_r||_{2q}$ and $\theta$ such that if the pinching condition
 \beqt\tag{$P_{K_{\theta}}$}
0\geqslant\lambda_1(M)\left( \int_MH_{r-1}dv_g\right)^2 -n\vvol(M)^{2}||H_r||_{2q}^2>-K_{\theta}
\eeqt
 is satisfied, then $M$ is diffeomorphic and $\theta$-quasi-isometric to $\Ss^n\left( \sqrt{\frac{n}{\lambda_1}}\right) $.
\end{thrm}
\begin{rem}
Note that if $q\geqslant\frac{n}{2r}$; then $ K_{\theta}$ does not depend on $||H_r||_{2q}$.
\end{rem}
\begin{rem}
The case $r=1$ was proved by Colbois and Grosjean \cite{CG}. In that case we do not need to assume that $H>0$.
\end{rem}
\section{Almost Einstein hypersurfaces}

We showed in \cite{Ro} that almost-Einstein hypersurfaces of $\R^{n+1}$ are close to round spheres. Namely,
\begin{thrm}[Roth \cite{Ro}]\label{thrm1}
Let $(M^n,g)$ be a compact, connected, oriented Riemannian manifold without boundary isometrically immersed in $\R^{n+1}$, $n\geqslant2$. Let $\theta\in]0,1[$. If $(M^n,g)$ is almost-Einstein, that is, $||\Ric-(n-1)kg||_{\infty}\leqslant\varepsilon$ for a positive constant $k$, with $\varepsilon$ small enough depending on $n$, $k$, $\Hinf$ and $\theta$, then $M$ is diffeomorphic and $\theta$-quasi-isometric to $\Ss^n\left(\sqrt{\frac{1}{k}}\right)$
\end{thrm}
By $\theta$-quasi-isometric, we understand that there exists a diffeomorphism $F$ from $M$ into $\Ss^n\left(\sqrt{\frac{1}{k}}\right)$ such that, for any $x\in M$ and for any unitary vector $u\in T_xM$, we have
$$\Big||d_xF(u)|^2-1\Big|\leqslant\theta.$$
This theorem is a corollary of our pinching result for the first eigenvalue of the Laplacian (Theorem A).\\
\indent
In this article, we consider almost-Einstein hypersurfaces of $\R^{n+1}$ in a weaker sense, namely for the $L^q$-norm, that is, $||\Ric-(n-1)kg||_q\leqslant\varepsilon$ for some positive constant $k$ and a sufficiently small $\varepsilon$. We prove that for some suitable constants $k$, such manifolds are close to round spheres. Precisely, we prove the following
\begin{thm}\label{thm1}
Let $(M^n,g)$ be a compact, connected, oriented Riemannian manifold without boundary isometrically immersed in $\R^{n+1}$. Let $q>\frac{n}{2}$, $r\in\{1,\cdots,n\}$ and if $r>1$, assume that $H_r>0$. Let $\theta\in]0,1[$, if $(M^n,g)$ satisfies $||\Ric-(n-1)k_{p,r}g||_q\leqslant\varepsilon$ for some sufficiently small $\varepsilon$ depending on $n$, $q$, $\Hinf$, $\vvol(M)$ and $\theta$, then $M$ is diffeomorphic and $\theta$-quasi-isometric to $\Ss^n\left( \sqrt{\frac{1}{k_{p,r}}}\right)$.  
\end{thm}
The constants $k_{p,r}$ in the theorem are defined from the higher order mean curvature (see Sect. \ref{prelim}).\\
\indent
After giving the proof of this theorem, we will give some applications to almost-umbilical hypersurfaces. \\ \\
{\it Proof:}
The proof is based on the above pinching result combined with a lower bound for the first eigenvalue of the Laplacian due to Aubry \cite{Au}. Assume that $||\Ric-(n-1)kg||_q\leqslant\varepsilon(n,q,k)$ for a positive constant $k$, $q>\frac{n}{2}$ and $\varepsilon$ small enough, then from Theorem 1.1 of \cite{Au}, we deduce that $\lambda_1(\Delta)$ satisfies
\beqt\label{lowerbound}
\lambda_1(\Delta)\geqslant nk(1-C_{\varepsilon}),
\eeqt
where $C_{\varepsilon}$ is an explicit constant such that $C_{\varepsilon}\lgra0$ when $\varepsilon\lgra0$.\\
\indent
Now, with the particular choice of $k=k_{p,r}$, we get:
 $$\lambda_1(M)\left( \int_MH_{r-1}dv_g\right)^2 -n\vvol(M)^{2}||H_r||_{2p}^2>-K_{\varepsilon}$$
for some constant $K_{\varepsilon}$ such that $K_{\varepsilon}\lgra0$ when $\varepsilon\lgra0$. \\
\indent
Let $\theta\in]0,1[$, we choose $\varepsilon(n,q,k,\theta)$ small enough such that $K_{\varepsilon}$ is small enough in Theorem \ref{thrm2} to obtain a diffeomorphism and $\theta$-quasi-isometry between $M$ and $\Ss^n\left( \sqrt{\frac{1}{k_{p,r}}}\right)$.$\finpreuve$ 

\begin{rem}
Note that in this Theorem, $\varepsilon$ depends on $\vvol(M)$ contrary to Theorem \ref{thrm1}. We can remove this dependence on the volume and replace it by a dependence on $||H_{r-1}||_1$. Indeed, using an upper bound on the volume under the $L^p$ condition on Ricci (result by Aubry \cite{Au}) we have that the volume is bounded from above by a constant depending on $n$, $q$, $\Hinf$ and $k_{r,p}$ that is in fact on $n$, $q$, $\Hinf$ and $||H_{r-1}||_1$ because of the definition of $k_{r,p}$ and \eqref{ineq}. On the other hand, by the classical extrisnic Sobolev inequality of Hoffman and Spruck \cite{HS}, we can bound the volume from below by a constant depending only on $n$ and $\Hinf$.
\end{rem}
Now, we will deduce from Theorem \ref{thm1} some Corollaries for almost-umbilical hypersurfaces of $\R^{n+1}$.
\section{Almost umbilical hypersurfaces}\label{secumb}
First, we give the following theorem, which a direct application of Theorem \ref{thrm1}.
\begin{thm}\label{thm2}
Let $(M^n,g)$ be a compact, connected, oriented Riemannian manifold without boundary isometrically immersed in $\R^{n+1}$. Let $\theta\in]0,1[$. If $(M^n,g)$ is almost-umbilical, that is, $||B-kg||_{\infty}\leqslant\varepsilon$ for a positive constant $k$, with $\varepsilon$ small enough depending on $n$, $k$ and $\theta$ then $M$ is diffeomorphic and $\theta$-quasi-isometric to $\Ss^n\left(\frac{1}{k}\right)$.
\end{thm}
\pf From the Gauss formula, we have
\beqt\label{gauss}
\Ric(Y,Y)=nH\left\langle B(Y),Y\right\rangle -\left\langle B(Y),B(Y)\right\rangle, 
\eeqt
for a tangent vector field $Y$. From (\ref{gauss}) and $||B-kg||_{\infty}\leqslant\varepsilon$, we deduce
\beQ
\Ric(Y,Y)&\geqslant& nk^2||Y||^2(1-\varepsilon)^2-k^2||Y||^2(1+\varepsilon)^2\\\\
&\geqslant& (n-1)k^2||Y||^2-\alpha_n(\varepsilon)||Y||^2,
\eeQ
where $\alpha_n$ is a positive function such that $\alpha_n(\varepsilon)\lgra0$ when $\varepsilon\lgra0$.\\
Similarly, we get
\beQ
\Ric(Y,Y)&\leqslant& (n-1)k^2||Y||^2+\alpha_n(\varepsilon)||Y||^2.
\eeQ
Finally, we have
$$||\Ric-(n-1)k^2g||_{\infty}\leqslant\alpha_n(\varepsilon),$$
which implies, by Theorem \ref{thrm1}, that for $\varepsilon$ small enough, $M$ is diffeomorphic and $\theta$-quasi-isometric to $\Ss^n\left( \frac{1}{k}\right)$.$\finpreuve$ \\
\indent
Now, from Theorem \ref{thm1}, it is possible to obtain, in some particular cases, comparable results for almost umbilical hypersurfaces in an $L^q$-sense. We recall that the umbilicity tensor is defined by
$$\tau=B-H\iid.$$
As we mentionned above, if $M$ is umbilical, {\it i.e.},$\tau=0$, and if $M$ is compact, then $M$ is a geodesic sphere. Here, we prove the following stability result for almost umbilical hypersurfaces.
\begin{thm}\label{thm3}
Let $(M^n,g)$ be a compact, connected, oriented Riemannian manifold without boundary isometrically immersed in $\R^{n+1}$. Let $q>\frac{n}{2}$, $r\in\{1,\cdots,n\}$ and if $r>1$, assume that $H_r>0$. For any $\theta\in]0,1[$, there exists two constants $\varepsilon_i(\theta,n,\Hinf,\vvol(M))$, $i=1,2$, such that if 
\be
\item $||\tau||_{2q}\leqslant\varepsilon_1$,
\item $||H^2-k_{p,r}||_q\leqslant\varepsilon_2$, for $p\geqslant4$ and $1\leqslant r\leqslant n$,
\ee 
then $M$ is diffeomorphic and $\theta$-quasi-isometric to $\Ss^n\left(\frac{1}{\sqrt{k_{p,r}}}\right)$
\end{thm} 
\begin{rem}
Note that for $r=1$, the result is due to Grosjean-Roth (see \cite{GR}), using a pinching result which involves only $H_1$ (see \cite{CG}).
\end{rem}
\noindent
{\it Proof.} Here again, from the Gauss formula, we have
$$\Ric=nHB-B^2.$$
From this, we deduce that 
\beQ
\Ric-(n-1)H^2g&=&nHB-B^2-(n-1)H^2g\\\\
&=&(n-2)H\tau-\tau^2,
\eeQ
which implies
\beQ
||\Ric-(n-1)kg||_q&\leqslant&||\Ric-(n-1)H^2g||_q+(n-1)\sqrt{n}||(H^2-k)||_q\\\\
&\leqslant&(n-2)\Hinf||\tau||_{2q}+||\tau||^2_{2q}+(n-1)\sqrt{n}||(H^2-k)||_q\\\\
&\leqslant&(n-2)\Hinf\varepsilon_1+\varepsilon_1^2+(n-1)\sqrt{n}\varepsilon_2
\eeQ
Now, we conclude by taking $\varepsilon_1$ and $\varepsilon_2$ small enough depending on $n$, $\Hinf$ and $\theta$ in order to apply Theorem \ref{thm1} and obtain the $\theta$-quasi-isometry. $\finpreuve$\\
Then, we deduce the following corollary which is to compare with Theorem \ref{thm2} .
\begin{cor}\label{cor1}
Let $(M^n,g)$ be a compact, connected, oriented Riemannian manifold without boundary isometrically immersed in $\R^{n+1}$. Let $\theta\in]0,1[$. If $(M^n,g)$ is almost-umbilic, that is, $||B-\sqrt{k_{p,r}}g||_{2q}\leqslant\varepsilon$, for $q>\frac{n}{2}$, with $\varepsilon$ small enough depending on $n$, $\Hinf$ and  $\theta$ then $M$ is diffeomorphic and $\theta$-quasi-isometric to $\Ss^n\left(\frac{1}{\sqrt{k_{p,r}}}\right)$.
\end{cor}
\pf A simple computation shows that 
$$||H^2-k_{p,r}||_{2q}\leqslant\alpha_1||B-\sqrt{k_{p,r}}g||_{2q},\ \text{and}$$
$$||\tau||_{2q}\leqslant\alpha_2||B-\sqrt{k_{p,r}}g||_{2q},$$
for two constants $\alpha_1$ and $\alpha_2$ depending on $n$ and $\Hinf$. Since we assume that $||B-\sqrt{k_{p,r}}g||_{2q}\leqslant\varepsilon$, we get 
\be
\item $||H^2-k_{p,r}||_{2q}\leqslant\alpha_1\varepsilon$,
\item $||\tau||_{2q}\leqslant\alpha_2\varepsilon$.
\ee
For $\varepsilon$ small enough, the assumptions of Theorem \ref{thm3} are satisfied and we can conclude that $M$ is diffeomorphic and quasi-isometric to $\Ss^n\left(\frac{1}{\sqrt{k_{p,r}}}\right)$.$\finpreuve$
\begin{rem}
We want to point out that this corollary is an improvement of Theorem \ref{thm2} only in some sense. Indeed, we improve the $L^{\infty}$-proximity to an $L^{2q}$-proximity, but this corollary is valid only for some special constants $k_{p,r}$ and not for any positive constant as in Theorem \ref{thm2}. 
\end{rem}

\end{document}